\documentclass{elsart3-1}

 \usepackage{graphics}
 \usepackage[T1]{fontenc}
 \usepackage{amsmath}
\usepackage[dvips]{epsfig}
\usepackage{amssymb}
\usepackage{dsfont}

\usepackage[english,francais]{babel}

\newtheorem{theorem}{Theorem}[section]

\newtheorem{e-proposition}[theorem]{Proposition}

\newtheorem{e-definition}[theorem]{Definition\rm}
\newtheorem{remark}{Remark}

\newtheorem{theoreme}{Th\'eor\`eme}[section]

\newtheorem{proposition}[theoreme]{Proposition}

\renewcommand{\theequation}{\arabic{equation}}
\newcommand{\Prob}{\mathds{P}}
\newcommand{\Nor}{\mathcal{N}}

\newcommand{\Inte}{\mathds{N}}
\newcommand{\Inter}{\mathds{Z}}
\newcommand{\LL}{\mathds{L}}
\newcommand{\Real}{\mathds{R}}
\newcommand{\Esp}{\mathds{E}}

\newcommand{\Proba}{\mathcal{P}}
\newcommand{\Loi}{\mathcal{L}}

\def\limitedistl{\renewcommand{\arraystretch}{0.5}
\begin{array}[t]{c}
\stackrel{{\Loi}}{\longrightarrow} \\
{\begin{array}{c}\scriptstyle  N\rightarrow\infty
\end{array}}
\end{array}\renewcommand{\arraystretch}{1}}

\def\limitedistlk{\renewcommand{\arraystretch}{0.5}
\begin{array}[t]{c}
\stackrel{{\Loi}}{\longrightarrow} \\
{\begin{array}{c}\scriptstyle  N\rightarrow\infty \\
\scriptstyle  |k'-k| \to \infty \end{array}}
\end{array}\renewcommand{\arraystretch}{1}}

\def\limiteprob{\renewcommand{\arraystretch}{0.5}
\begin{array}[t]{c}
\stackrel{{\Proba}}{\longrightarrow} \\
{\scriptstyle N\rightarrow\infty}
\end{array}\renewcommand{\arraystretch}{1}}

\def\limiteNN{\renewcommand{\arraystretch}{0.5}
\begin{array}[t]{c}
\stackrel{}{\longrightarrow} \\
{\scriptstyle N\rightarrow\infty}
\end{array}\renewcommand{\arraystretch}{1}}

\def\og{\leavevmode\raise.3ex\hbox{$\scriptscriptstyle\langle\!\langle$~}}
\def\fg{\leavevmode\raise.3ex\hbox{~$\!\scriptscriptstyle\,\rangle\!\rangle$}}

\journal{the Acad\'emie des sciences}
\begin{document}
\centerline{}
\begin{frontmatter}


\selectlanguage{english}
\title{Detecting abrupt changes of the long-range dependence or the
self-similarity of a Gaussian process}


\selectlanguage{english}
\author[authorlabel1]{Jean-Marc Bardet},
\ead{jean-marc.bardet@univ-paris1.fr}
\author[authorlabel1]{Imen Kammoun}
\ead{imen.kammoun@univ-paris1.fr}

\address[authorlabel1]{Universit\'e Paris 1, SAMOS-MATISSE-CES, 90 rue de Tolbiac, 75013 Paris,
France}


\medskip
\begin{center}
{\small Received *****; accepted after revision +++++\\
Presented by \pounds\pounds\pounds\pounds\pounds}
\end{center}
\vspace{-0.5cm}
\begin{abstract}
\selectlanguage{english}
In this paper, an estimator of $m$ instants ($m$ is known) of abrupt
changes of the parameter of long-range dependence or self-similarity
is proved to satisfy a limit theorem with an explicit convergence
rate for a sample of a Gaussian process. In each estimated zone
where the parameter is supposed not to change, a central limit
theorem is established for the parameter's (of long-range
dependence, self-similarity) estimator and a goodness-of-fit test is
also built. {\it To cite this article: J.M. Bardet, I. Kammoun, C.
R. Acad. Sci. Paris, Ser. I 340 (2007).}

\vskip 0.5\baselineskip

\selectlanguage{francais}
\noindent{\bf R\'esum\'e} \vskip 0.5\baselineskip \noindent {\bf
D\'etection de ruptures du param\`etre de longue m\'emoire,
d'autosimilarit\'e pour des processus gaussiens.} Dans ce papier,
pour une trajectoire d'un processus gaussien, un estimateur des $m$
points de ruptures ($m$ est suppos\'e connu) du param\`etre de
longue m\'emoire ou d'autosimilarit\'e est construit et on montre
qu'il v\'erifie un th\'eor\`eme limite avec une vitesse de
convergence explicite. Dans chaque zone (estim\'ee) o\`u ce
param\`etre est constant, un estimateur de ce param\`etre vérifie un
th\'eor\`eme limite centrale et un test d'ajustement est également
mis en place. {\it Pour citer cet article~: J.M. Bardet, I. Kammoun,
C. R. Acad. Sci. Paris, Ser. I 340 (2007).}

\end{abstract}
\end{frontmatter}
\vspace{-1.2cm}
\selectlanguage{francais}
\section*{Version fran\c{c}aise abr\'eg\'ee}
Le probl\`eme de d\'etection de points de ruptures moyennant la
minimisation d'une fonction de contraste donn\'ee a \'et\'e
\'etudi\'e depuis le milieu des ann\'ees 1990 dans le cadre de
processus \`a longue m\'emoire (voir par exemple \cite{gls},
\cite{kl}, \cite{lavielle0}, \cite{lavielle1} et \cite{lavielle2}).
De ces approches, certaines ont \'et\'e associ\'ees \`a un cadre
param\'etrique, telle que la d\'etection de rupture selon la moyenne
et/ou la variance, d'autres trait\'ees dans un cadre non
param\'etrique (comme la d\'etection de ruptures selon la
distribution ou le spectre). Dans la litt\'erature, diff\'erents
auteurs ont \'egalement propos\'e des statistiques de test de
l'hypoth\`ese que le param\`etre est inchang\'e contre le fait que
le param\`etre de longue m\'emoire varie en fonction du temps (voir
par exemple \cite{aya}, \cite{bt}, \cite{hor}). \`A notre
connaissance, le cadre semi-param\'etrique de d\'etection de
changements en longue m\'emoire ou en autosimilarit\'e n'a été
traité que dans \cite{lavielle1} à partir d'une technique basée sur
le périodogramme.
\\
Notre approche est fondée sur l'analyse par ondelettes, ce qui présente plusieurs avantages:
c'est une technique non-paramétrique applicable pour des processus très généraux, robuste aux tendances polynomiales et, au moins dans le cadre gaussien,
s'accompagnant de tests d'adéquation de type $\chi^2$ simples et intéressants à utiliser. Ainsi, un estimateur des $m$ points de ruptures ($m \in
\Inte^*$, supposé connu) de la longue m\'emoire ou d'auto-similarit\'e est
con\c{c}u pour un \'echantillon de processus gaussien en se basant
sur l'analyse par ondelettes, ce qui permet ensuite de mettre en place des tests d'adéquation. Pour ce type de processus, cette
m\'ethode a \'et\'e propos\'ee pour la premi\`ere fois dans
\cite{flan2}, puis d\'evelopp\'ee par exemple dans \cite{avf}. La convergence des estimateurs bas\'es sur les
ondelettes a \'et\'e \'etudi\'ee dans le cas du mouvement brownien
fractionnaire (FBM) dans \cite{JM1}, et dans un cadre semi-param\'etrique g\'en\'eral de
processus  gaussiens stationnaires \`a longue
m\'emoire par \cite{mrt} et
\cite{bardet}. \\
Ici le principe de l'estimation du param\`etre de
longue m\'emoire ou d'auto-similarit\'e est le suivant: dans chaque
zone o\`u il n'y a pas de changement, ce param\`etre peut \^etre
estim\'e \`a partir d'une log-log r\'egression de la variance des
coefficients d'ondelettes sur plusieurs \'echelles choisies (voir
(\ref{TLClogS})). Une fonction de contraste d\'efinie par la somme
des carr\'es des distances entre ces points et les droites
d'ajustement, dans les $m+1$ zones possibles d\'etect\'ees, est
minimis\'ee (voir (\ref{contU})), donnant un estimateur des points
de ruptures (voir (\ref{paramestim})). Sous certaines hypoth\`eses
g\'en\'erales, on montre qu'il v\'erifie un th\'eor\`eme limite avec
une vitesse de convergence explicite (voir Theorem \ref{convprob}).
Dans chacune des zones d\'etect\'ees, les param\`etres de longue
m\'emoire, (ou d'auto-similarit\'e) peuvent \^etre estim\'es, tout
d'abord avec la regression des moindres carr\'es ordinaires (OLS),
puis par une regression des moindres carr\'es
pseudo-g\'en\'eralis\'es (FGLS). Un th\'eor\`eme de la limite
centrale est \'etabli pour chacun des deux estimateurs (voir Theorem
\ref{convparam} et Proposition \ref{FGLS} ci-dessous) et des
intervalles de confiance peuvent \^etre calcul\'es. L'estimateur
FGLS offre deux avantages: d'une part, sa variance asymptotique est
plus petite que celle de l'estimateur OLS, et d'autre part, il
permet la construction d'un test d'ajustement tr\`es simple bas\'e sur le carr\'e des distances entre les points
(d'abscisse, le logarithme d'une \'echelle choisie et d'ordonn\'ee,
le logarithme de la variance empirique des coefficients d'ondelettes
pour cette \'echelle) et les droites de r\'egression
pseudo-g\'en\'eralis\'ee correspondantes (voir (\ref{test})). La
convergence vers une distribution du Chi-deux de
ce test est \'etablie dans le Theorem \ref{T}. \\
Deux cas particuliers de processus gaussiens sont ensuite
\'etudi\'es dans la section \ref{simu}. En premier lieu, on
s'int\'eresse aux s\'eries chronologiques stationnaires longue
m\'emoire avec un param\`etre de Hurst constant par morceaux. On se
place dans un cadre semi-param\'etrique contenant par exemple les
FGN et les processus FARIMA (voir Figure \ref{Figure}). En second
lieu, le cas d'un processus \`a accroissements stationnaires et
autosimilaire par morceaux est trait\'e, ce qui revient \`a
consid\'erer des successions de FBM ayant des exposants de Hurst
distincts (voir Figure \ref{Figure}).
Pour ces deux exemples, les vitesses de convergence explicites des
diff\'erents estimateurs et tests sont donn\'ees et des simulations
montrent leurs qualit\'es (voir Table \ref{Table}). D'autres
simulations, preuves des théorèmes ainsi qu'un exemple plus général
de détection de ruptures dans le cadre de processus gaussien
localement fractionnaire sont détaillés dans \cite{bardet3}.
\vspace{-0.2cm}

\selectlanguage{english}
\vspace{-0.67cm}
\section{Assumptions and main results}
Let $(X_t)_{t \in \Inte}$ be a Gaussian process and assume
that $\big (X_0,X_1,\ldots,X_{N}\big )$ is known.
In the sequel, $X$ will be a piecewise stationary long memory time
series or a piecewise self-similar time series having stationary
increments. Consider $\psi:\Real \rightarrow \Real$ a function
called "the mother wavelet". For $(a,b)\in \Real_+^*\times \Real$,
the wavelet coefficient of $X$ for the scale $a$ and the shift $b$
is $
d_{X}(a,b):=\frac{1}{\sqrt{a}}\int_{\Real}\psi(\frac{t-b}{a})X(t)dt$.
When only a discretized path of $X$ is available, approximations
$e_{X}(a,b)$ are only computable: , $$e_{X}(a,b):=\frac {1}{\sqrt a}
\, \sum_{p=1}^{N} \psi\big (\frac {p-b} {a_N} \big ) X_{p
}~~~\mbox{for $(a,b)\in \Real_+^*\times \Inte.$}
$$
~\\
Assume that there exist $m\in \Inte$ (the number of abrupt changes)
and\vspace{2mm}
\begin{itemize}
\item [$\bullet$] $0=\tau_0^*<\tau_1^*<\ldots<\tau_{m}^*<\tau_{m+1}^*=1$ (unknown
parameters);
\item [$\bullet$] two families $(\alpha_j^*)_{0\leq j \leq m}\in
\Real^{m+1}$ and $(\beta_j^*)_{0\leq j \leq m}\in (0,\infty)^{m+1}$
(unknown parameters);
\item [$\bullet$] a sequence of "scales" $(a_n)_{n \in
\Inte}\in \Real^\Inte$ (chosen) satisfying $a_n \geq a_{min}$, with
$a_{min}>0$,\\
\end{itemize}
\vspace{-2mm} $\mbox{such that for } j=0,1,\ldots,m \mbox{ and } k
\in D_N^*(j) \subset \big [[N\tau_j^*],[N\tau_{j+1}^*]\big ],$
$$\Esp \big
[e_{X}^2(a_N,k)\big ] \sim \beta_j^* \big ( a_N\big )^{\alpha_j^*}
\mbox{ when } N\to \infty.
$$ A piecewise sample variance can be
the appropriated estimator of such power law. Thus, define
$$
S^{k'}_{k}(a_N):= \frac {a_N} {k'-k} \, \sum
_{p=[k/a_N]}^{[k'/a_N]-1}e_{X}^2(a_N,a_N\,p)~~~\mbox{for $0\leq k <
k'\leq N$}. $$ Now set $0<r_1<\ldots<r_\ell$ with $\ell \in
\Inte^*$, and assume that a multidimensional central limit theorem
can be established for $\Big(\log\big(S^{k'}_{k}(r_i \, a_N)\big
)\Big)_{1 \leq i \leq \ell}$, when $[N\tau_j^*] \leq k < k' \leq
[N\tau_{j+1}^*]$, {\it i.e.}\vspace{-2mm}
\begin{eqnarray}\label{TLClogS}
\hspace{-1cm} \sqrt {\frac {k'-k} {a_N }} \Big (\log \big
(S^{k'}_{k}(r_i \, a_N)\big )-\log (\beta_j^*)-\alpha_j^* \log \big
(r_i \, a_N\big ) \Big )_{1 \leq i \leq \ell} \limitedistlk
\Nor\big(0,\Gamma^{(j)}(\alpha_j^*,r_1,\ldots,r_\ell)\big),\vspace{-4mm}
\end{eqnarray}
with $\Gamma^{(j)}(\alpha_j^*,r_1,\ldots,r_\ell)=\big (
\gamma_{pq}^{(j)}\big ) _{1 \leq p,q \leq \ell}$ a $(\ell \times
\ell)$ matrix not depending on $N$ such that $\alpha \mapsto
\Gamma^{(j)}(\alpha,r_1,\ldots,r_\ell)$ is a continuous function and
a positive matrix for all $\alpha$. Define a contrast
function\vspace{-2mm}
\begin{eqnarray}\label{contU}
\hspace{-5mm}U_N\big((\alpha_j)_{0\leq j\leq m},\,(\beta_j)_{0\leq
j\leq m},\, (k_j)_{1\leq j\leq
m}\big)=\sum_{j=0}^{m}\sum_{i=1}^{\ell} \Big(
\log\big(S^{k_{j+1}}_{k_j}(r_i \,a_N)\big)- \big(\alpha_j\log (r_i
\,a_N)+\log \beta_j\big )\Big)^2\vspace{-1mm}
\end{eqnarray}
with $(\alpha_j)_{0\leq j\leq m}\in A^{m+1}\subset \Real^{m+1}$,
$(\beta_j)_{0\leq j\leq m} \in B^{m+1} \subset(0,\infty)^{m+1}$,
$0=k_0<k_1<\ldots<k_m<k_{m+1}=N, (k_j)_{1\leq j\leq m} \in K_{m}(N)
\subset \Inte^{m}$. The vector of estimated parameters $\widehat
\alpha_j,~\widehat \beta_j$ and $\widehat k_j$ (and therefore
$\widehat \tau_j$) is the vector which minimizes this contrast
function in $A^{m+1}\times B^{m+1}\times K_{m}(N)$, {\it i.e.},
\begin{eqnarray}\label{paramestim}
\big ((\widehat\alpha_j)_{0\leq j\leq m},\,(\widehat \beta_j)_{0\leq
j\leq m},\, (\widehat k_j)_{1\leq j\leq m}\big )&:=&\text{Argmin}
\Big\{U_N\big((\alpha_j)_{0\leq j\leq m},\,(\beta_j)_{0\leq j\leq
m},\, (k_j)_{1\leq j\leq
m}\big)\Big\}\\
\widehat \tau_j&:=&\widehat k_j/N~~\mbox{for}~~1\leq j\leq m.
\end{eqnarray}
For a given $(k_j)_{1\leq j\leq m}$, it is obvious that
$(\widehat\alpha_j)_{0\leq j\leq m}$ and $( \log
\widehat\beta_j)_{0\leq j\leq m}$ are obtained from a log-log
regression of $\big (S^{k_{j+1}}_{k_j}(r_i \, a_N)\big )_i$ onto
$\big (r_i \, a_N\big)_i$, {\it i.e.}\vspace{-3mm}
$$
\left (\hspace{-1mm} \begin{array}{c} \widehat \alpha_j \\
\log \widehat \beta_j \end{array} \hspace{-1mm}\right ) =\big (
L_{1}' \cdot L_{1})^{-1} L_{1}'\cdot Y_{k_j}^{k_{j+1}}
~~\mbox{with}~~
Y_{k_j}^{k_{j+1}}\hspace{-1mm}:=\big(\log\big(S^{k_{j+1}}_{k_j}(r_i\cdot
a_N)\big)\big)_{1\leq i\leq \ell}~~,~~ L_{a_N}:=\left  (
\begin{array}{cc} \log (r_1\, a_N) & 1\\ \vdots \vspace{-1mm}&
\vdots\vspace{-1mm} \\\log (r_\ell \, a_N) & 1
\end{array} \right).
$$
Therefore $(\widehat k_j)_{1\leq j\leq m}=\text{Argmin}
\Big\{U_N\big((\widehat \alpha_j)_{0\leq j\leq m},\,(\widehat
\beta_j)_{0\leq j\leq m},\, (k_j)_{1\leq j\leq
m}\big),~~(k_j)_{1\leq j\leq m} \in K_{m}(N)\Big \}$.\\
\begin{remark}
In this paper, $m$ is supposed to be known. However, if $m$ is
unknown, as in \cite{lavielle0} or \cite{lavielle1}, a penalized
contrast $\tilde U_{m,N}=U_N+\beta_N \times  m$ (with $\beta_N$ an
appropriated sequence converging to $0$) can be used instead of
$U_N$, and by adding a minimization in $m$, an estimator $\widehat
m$ of $m$ could be also deduced.\\\vspace{-3mm}
\end{remark}
In this paper, parameters $(\alpha_j^*)$ are supposed to satisfied
abrupt changes:\\
{\bf Assumption C :} Parameters $(\alpha_j^*)$ are such that $
|\alpha_{j+1}^*-\alpha_j^*|\neq 0~~\mbox{for
all}~~j=0,1,\ldots,m-1$.
\begin{theorem}\label{convprob} Define $
\underline{\tau}^*:=(\tau^*_1,\ldots,\tau^*_m),~~\underline{\widehat\tau}:=(\widehat
\tau_1,\ldots,\widehat \tau_m)$ and $\|\underline{\tau}\|_m:=\max
\big (|\tau_1|,\ldots,|\tau_m|\big )$. Let $\ell \in \Inte \setminus
\{0,1,2\}$. If Assumption C and relation (\ref{TLClogS}) holds with
$(\alpha_j^*)_{0\leq j \leq m}$ such that $\alpha_j^* \in [a\, , \,
a']$ and $a<a'$ for all $j=0,\ldots,m$, then if $\displaystyle~~
a_N^{1+2(a'-a)} N ^{-1} \limiteNN 0$, for all $(v_n)_n$ satisfying
$\displaystyle~~ v_N \cdot a_N^{1+2(a'-a)} N^{-1} \limiteNN
0$,\vspace{-3mm}
\begin{eqnarray}\label{conv_tau}
\hspace{2cm} \Prob \Big (v_N
\|\underline{\tau}^*-\underline{\widehat \tau}\|_m \geq \eta \Big )
\limiteNN 0~~~\mbox{for all $\eta>0$}.\vspace{-5mm}
\end{eqnarray}
\end{theorem}
\begin{remark}
The proof of this result is provided in \cite{bardet3}. Unfortunately, the rate of convergence of $\|\underline{\tau}^*-\underline{\widehat \tau}\|_m$ is only $v_N=N^\alpha$ with $0<\alpha<1$ and not $N$ as, for instance, in
\cite{lavielle0} and \cite{lavielle1}. However the context is not the same: in these papers, the contrast is directly computed from $N$ values of $(X)_i$ which do not change following $\widehat \tau$. Here, the contrast is computed from only $(m+1)\ell$ values of $S_N$ which change  following $\widehat \tau$. The rate of convergence $N$ can not be reached in such a context (simulations show also this property). This is certainly a drawback of your method, which hopefully does not change the rate of convergence of parameters $(\alpha_j)$ and $(\beta_j)$.
\end{remark}
For
$j=0,1,\ldots,m$, the log-log regression of $\big(S_{\widehat
k_j}^{\widehat k_{j+1}}(r_i a_N)\big)_{1\leq i\leq\ell}$ onto $(r_i
a_N)_{1\leq i\leq\ell}$ provides estimators of $\alpha_j^*$ and
$\beta_j^*$. However, if $ \tau_j$ converges to $\tau_j^*$,
$\widehat k_j=N\cdot \widehat \tau_j$ does not converge to $k_j^*$,
and therefore $\Prob\big([\widehat k_j,\widehat
k_{j+1}]\subset[k^*_j,k^*_{j+1}]\big)$ does not tend to $1$. So,
define $\tilde k_j$ and $\tilde k_{j}'$ such that $ \tilde
k_j=\widehat k_j +\frac {N} {v_N}$ and $\tilde k_{j}'=\widehat
k_{j+1}-\frac {N} {v_N}$. From (\ref{conv_tau}) with $\eta=1/2$,
$\Prob\big([\tilde k_j,\tilde k_{j}']\subset[k^*_j,k^*_{j+1}]\big)
\limiteNN 1$. Then,\\
\vspace{-3mm}
\begin{theorem}\label{convparam}
Let $\Theta^*_{j}:=\Big(
\begin{array}{c}
\alpha^*_j \\
\log \beta^*_j\\
\end{array} \Big)$ and $\tilde\Theta_{j}:=(L_{1}^{'} L_{1})^{-1} L_{1}^{'}
Y_{\tilde k_j}^{\tilde k_{j}'}=\Big(
\begin{array}{c}
\tilde \alpha_j \\
\log \tilde \beta_j\\
\end{array} \Big )$. Under the same assumptions as in Theorem \ref{convprob}, for
$j=0,\ldots,m$, with
$\Sigma^{(j)}(\alpha_j^*,r_1,\ldots,r_\ell):=(L_{1}^{'} L_{1})^{-1}
L_{1}^{'} \Gamma^{(j)}(\alpha_j^*,r_1,\ldots,r_\ell) L_{1}
(L_{1}^{'}  L_{1})^{-1}$,\vspace{-2mm}
\begin{eqnarray}\label{conv_theta}
\hspace{2cm} \sqrt {\frac {N \big(\tau^*_{j+1}-\tau^*_j\big)}{a_N}}
\Big(\tilde\Theta_{j}-\Theta^*_{j}\Big)\limitedistl
\Nor\big(0,\Sigma^{(j)}(\alpha_j^*,r_1,\ldots,r_\ell)\big)\vspace{-3mm}
\end{eqnarray}
\end{theorem}
A second estimator of $\Theta^*_{j}$ can be obtained from feasible
generalized least squares (FGLS) estimation. Indeed, the asymptotic
covariance matrix $\Gamma^{(j)}(\alpha_j^*,r_1,\ldots,r_\ell)$ can
be estimated by the matrix $\tilde \Gamma^{(j)}:=
\Gamma^{(j)}(\tilde \alpha_j,r_1,\ldots,r_\ell)$ and $ \tilde
\Gamma^{(j)}\limiteprob \Gamma^{(j)}(\alpha_j^*,r_1,\ldots,r_\ell)$.
Then, the FGLS estimator ${\overline {\Theta}}_j$ of $\Theta^*_{j}$
is defined from the minimization among all $\Theta $ of the following
squared distance,
$$
\parallel Y_{\tilde k_j}^{\tilde
k_{j}'}-L_{a_N}\cdot \Theta
\parallel^2_{\tilde \Gamma^{(j)}}=\big(Y_{\tilde k_j}^{\tilde
k_{j}'}-L_{a_N} \Theta\big)'\cdot  \big(\tilde
\Gamma^{(j)}\big)^{-1}\cdot \big(Y_{\tilde k_j}^{\tilde
k_{j}'}-L_{a_N}  \Theta\big).
$$
and therefore define $ {\overline {\Theta}}_j:=\big(L_{1}'
\big(\tilde \Gamma^{(j)}\big)^{-1} L_{1} \big)^{-1} L_{1}'
 \big(\tilde \Gamma^{(j)}\big)^{-1} Y_{\tilde
k_j}^{\tilde k_{j}'}$.\\
\vspace{-2mm}
\begin{proposition}\label{FGLS}
Under the same assumptions as in Theorem \ref{convparam}, for
$j=0,\ldots,m$\vspace{-0.5mm}
\begin{eqnarray} \label{estim2}
\hspace{2cm} \sqrt {\frac { N \big(\tau^*_{j+1}-\tau^*_j\big)}{a_N}}
\Big({\overline {\Theta}}_j-\Theta^*_{j}\Big)\limitedistl
\Nor\big(0,M^{(j)}(\alpha_j^*,r_1,\ldots,r_\ell)\big)\vspace{-3mm}
\end{eqnarray}
with $M^{(j)}(\alpha_j^*,r_1,\ldots,r_\ell):=\big(L_{1}^{'}
\big(\Gamma^{(j)}(\alpha_j^*,r_1,\ldots,r_\ell)\big)^{-1}
L_{1}\big)^{-1}\leq \Sigma^{(j)}(\alpha_j^*,r_1,\ldots,r_\ell)$
(for the order's relation between positive symmetric matrix).\\
\end{proposition}
Therefore ${\overline {\alpha}}_j$ is more accurate than $\tilde
\alpha_j$ for estimating $\alpha_j^*$ when $N$ is large enough. For
$j=0,\ldots,m$, let $T^{(j)}$ be the FGLS distance between points
$\Big(\log(r_i\, a_N),\log\big(S^{\tilde k_{j}'}_{\tilde
k_j}\big)\Big)_{1\leq i\leq\ell}$ and the FGLS regression line. The
following theorem describes  the asymptotic behavior of a
goodness-of-fit test on each segment $[\tilde k_j,\tilde k_{j}'[$:
\begin{theorem}\label{T}
Under the same assumptions as in Theorem \ref{convprob}, for
$j=0,\ldots,m$\vspace{-2mm}
\begin{eqnarray}\label{test}
T^{(j)}= \frac {N \big(\tau^*_{j+1}-\tau^*_j\big)}{a_N}
\parallel Y_{\tilde k_j}^{\tilde
k_{j}'}-L_{a_N}  {\overline {\Theta}}_j\parallel^2_{\tilde
\Gamma^{(j)}}\limitedistl \chi^2(\ell-2).\vspace{-5mm}
\end{eqnarray}
\end{theorem}
\vspace{-1cm}
\section{Applications}\label{simu}
\vspace{-0.2cm}
\subsection{Piecewise long memory Gaussian processes}
\vspace{-0.2cm} Assume that the process $X=(X_t)_{t\in \Inte}$ is a
Gaussian piecewise long-range dependent (LRD) process, {\it i.e.}
there exists $(D_j^*)_{0\leq j\leq m} \in (0,1)^{m+1}$ and for all
$j=0,\ldots,m$ and $k \in \big
\{[N\tau_j^*],[N\tau_j^*]+1,\ldots,[N\tau_{j+1}^*]-1\big \}$,
$X_k=X^{(j)}_{k-[N\tau_{j}^*]}$, where $X^{(j)}=(X^{(j)}_t)_{t\in
\Inte}$ satisfies the following Assumption LRD$(D_j^*)$.\vspace{2mm}\\
{\bf Assumption LRD$(D)$:} $Y$ is a centered stationary Gaussian
process with spectral density $f$ such that $f(\lambda)= |\lambda
|^{-D} \cdot f^*(\lambda)~~\mbox{for all}~~ \lambda \in
[-\pi,\pi]\setminus \{0\}$ with $f^*(0)>0$ and with $C_2>0$,
$|f^*(\lambda)-f^*(0)|\leq C_{2} \cdot |\lambda|^{2}~~\mbox{for
all}~\lambda \in [-\pi,\pi]$.\vspace{2mm} \\
Following \cite{bardet}, if the mother wavelet is supposed to be
included in a Sobolev ball, then \\ \vspace{-1mm}
\begin{cor}\label{cor_LRD}
Let $X$ be a Gaussian piecewise LRD process defined as above and
$\psi:~\Real \mapsto \Real$ be $[0,1]$-supported with
$\psi(0)=\psi(1)=0$ and $ \int_0^1 \psi(t)\,dt=0$ and such that
there exists sequence $(\psi_\ell)_{\ell \in \Inter}$ satisfying
$\psi(\lambda)=\sum_{\ell \in \Inter}\psi_\ell e^{2\pi i\ell
\lambda} \in\LL^2([0,1])$ and $\sum_{\ell\in\Inter}(1+|\ell|)^{5/2}
|\psi_\ell|<\infty$.  Under Assumption C, for all $0<\kappa< 2/15$,
if $a_N=N^{\kappa+1/5}$ and $v_N=N^{2/5-3\kappa}$ then
(\ref{conv_tau}), (\ref{conv_theta}), (\ref{estim2}) and
(\ref{test}) hold.\\
\end{cor}
Thus, the rate of convergence of $\underline{\widehat \tau}$ to
$\underline{\tau}^*$ (in probability) is $N^{2/5-3\kappa}$ for
$0<\kappa$ arbitrary small. Estimators $\tilde D_j$ and
${\overline D_j}$ converge to the parameters $D_j^*$ following a
central limit theorem with a rate of convergence $N^{2/5-\kappa/2}$
for $0<\kappa$. Convincing results of
simulations can be observed in Table \ref{Table} and Figure
\ref{Figure}.\\
\vspace{-0.3cm}
\begin{table}[!h]
\renewcommand{\arraystretch}{1.2}
\begin{center}
\begin{tabular}{ccccc}
\begin{tabular} {c}
\vspace{0.3mm}
\\
\vspace{0.5mm} \vspace{0.5mm}
\\
\hline
\begin{tabular} {|c|}
~\emph{Estim.}~
\\
\\
\hline ~$\widehat \sigma_{\mbox{\emph{Estim.}}}$~\\
\hline ~$\sqrt{MSE}$~
\end{tabular}
\\
\hline
\end{tabular}
&&
\begin{tabular} {|c|}
\hline
$N=20000$\\
\hline \hline \multicolumn{1}{|c|}{$\tau_1$:0.75  $D_0$:0.2  $D_1$:0.8}\\
\hline \hline
\begin{tabular} {c|c|c}
$\widehat \tau_1$ & $\tilde D_0$ & $\tilde D_1$\\
~0.7540~ & ~0.1902~ & ~0.7926~ \\
\hline
~0.0215~ & ~0.0489~  & ~0.0761~   \\
\hline
~0.0218~ & ~0.0499~ & ~0.0764~ \\
\end{tabular}
\\
\hline
\end{tabular}
&&
\begin{tabular} {|c||c|}
\hline
$N=5000$ & $N=10000$\\
\hline \hline \multicolumn{2}{|c|}{$\tau_1$:0.3 ~~~~ $\tau_2$:0.78
~~~~ $H_0$:0.6 ~~~~ $H_1$:0.8 ~~~~ $H_2$:0.5}\\
\hline \hline
\begin{tabular} {c|c|c|c|c}
$\widehat \tau_1$ & $\widehat \tau_2$ & $\tilde H_0$ & $\tilde
H_1$ & $\tilde H_2$ \\
~0.3465~ & ~0.7942~ & ~0.5578~ & ~0.7272~ & ~0.4395~ \\
\hline
~0.1212~  & ~0.1322~ & ~0.0595~ & ~0.0837~ & ~0.0643~ \\
\hline
~0.1298~ & ~0.1330~ & ~0.0730~ & ~0.1110~ & ~0.0883~\\
\end{tabular}
&
\begin{tabular} {c|c|c|c|c}
$\widehat \tau_1$ & $\widehat \tau_2$ & $\tilde H_0$ & $\tilde
H_1$ & $\tilde H_2$ \\
~0.3086~  & ~0.7669~ & ~0.5597~ & ~0.7633~ & ~0.4993~ \\
\hline
~0.0893~  & ~0.0675~ & ~0.0449~ & ~0.0813~ & ~0.0780~ \\
\hline
~0.0897~ & ~0.0687~ & ~0.0604~ & ~0.0892~ & ~0.0780~\\
\end{tabular}\\
\hline
\end{tabular}
\end{tabular}
\end{center}
\caption{Left: Estimation of $\tau_1$, $D_0$ and $D_1$ in the case
of piecewise FARIMA(0,$d_j$,0) ($d_1=0.1$ and $d_2=0.4$) with one
change point when $N=20000$ (50 realizations). Right: Estimation of
$\tau_1$, $\tau_2$, $H_0$, $H_1$ and $H_2$ in the case of piecewise
FBM with two change points when $N=5000$ and $N=10000$ (50
realizations)}\label{Table}
\end{table}
\vspace{-0.6cm}
\subsection{Piecewise fractional Brownian motions}
\begin{figure}[h!]
 \begin{minipage}[b]{.45\linewidth}
  \raggedleft\includegraphics[width=7.5 cm,height=3.6 cm]{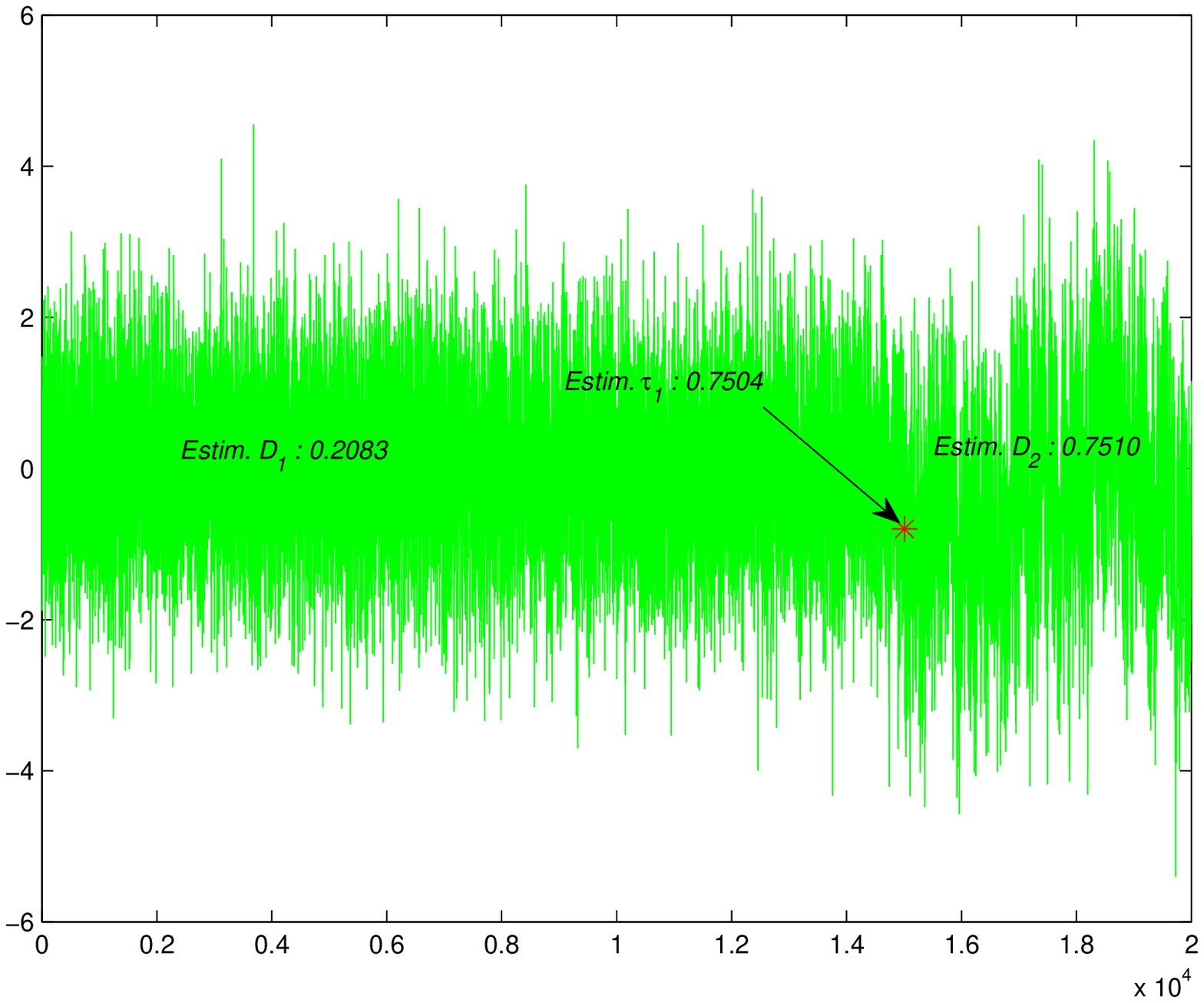}
 \end{minipage} \hfill
 \begin{minipage}[b]{.5\linewidth}
  \raggedright\includegraphics[width=7 cm,height=3.4 cm]{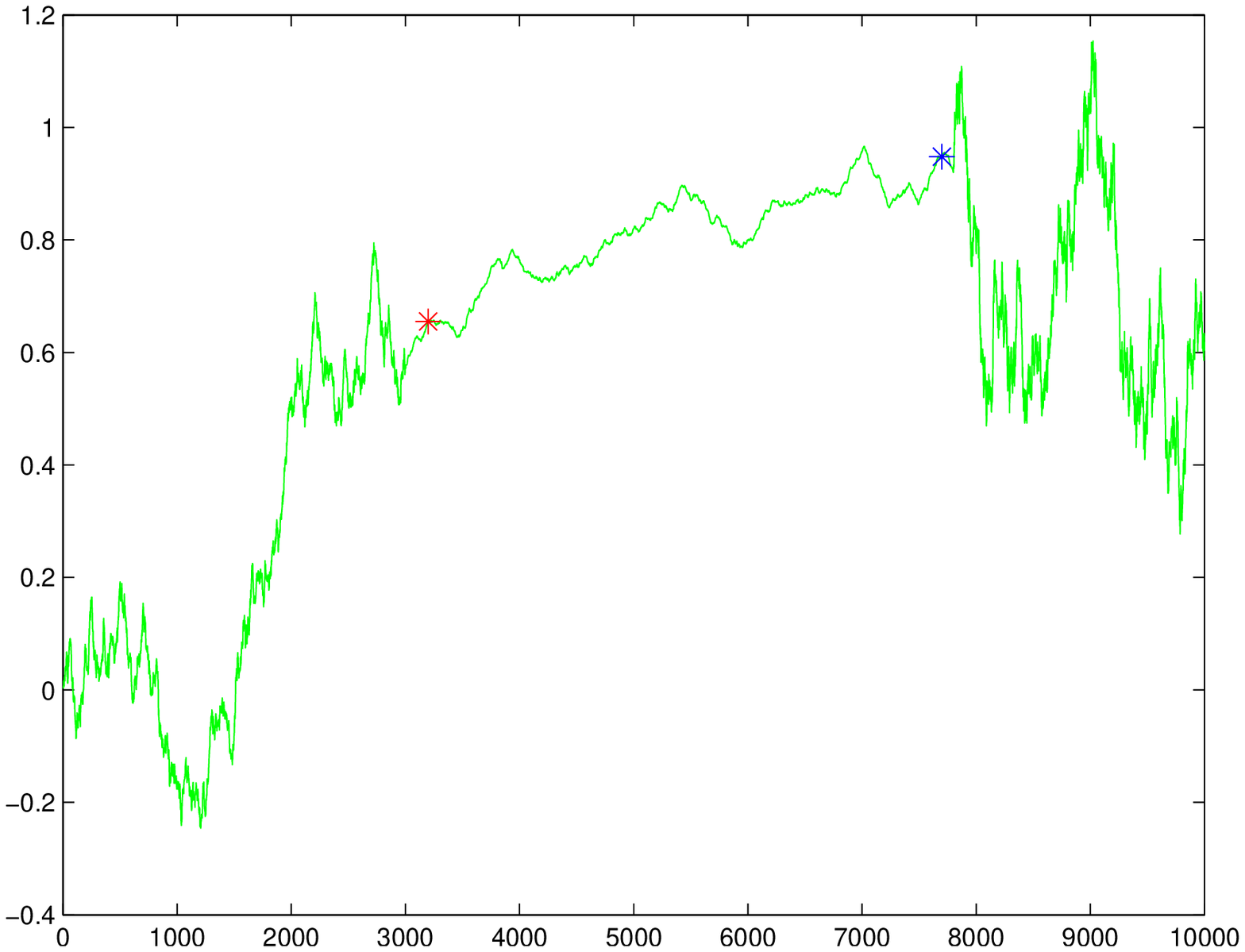}
  \end{minipage}\vspace{-0.3cm}
 \caption{Left: Piecewise FARIMA(0,$d_j$,0) (with $d_0:0.1$ ($D_0:0.2$), $d_1:0.4$ ($D_1:0.8$) and
  $\tau_1:0.75$). Right: Piecewice FBM($H_j$) ($\tau_1:0.3$, $\tau_2:0.78$, $H_0:0.6$, $H_1:0.8$ and $H_2:0.5$),
($\widehat\tau_1:0.32$, $\widehat\tau_2:0.77$, $\tilde H_0:0.5608$,
$\tilde H_1:0.7814$ and $\tilde H_2:0.4751$). \label{Figure}}
\end{figure}

\vspace{-0.2cm} Now, $X$ will be called a piecewise fractional
Brownian motion if there exist two families of parameters
$(H^*_j)_{0\leq j\leq m}\in (0,1)^{m+1}$ and $(\sigma^{*2}_j)_{0\leq
j\leq m}\in (0,\infty)^{m+1}$ such that for all $j=0,\ldots,m$ and
$t \in \big [[N\tau_j^*],[N\tau_j^*]+1,\ldots,[N\tau_{j+1}^*]-1\big
]$, $X_t=X^{(j)}_{t-[N\tau_{j}^*]}$, where
$X^{(j)}=(X^{(j)}_t)_{t\in \Real}$ is a FBM with parameters $H^*_j$
and $\sigma^{*2}_j$. Following the results of \cite{JM1}, one
obtains,\\ \vspace{-1mm}
\begin{cor}\label{cor_FBM}
Let $X$ be a piecewise FBM and $\psi:\Real \to \Real$ be a piecewise
continuous and left (or right)-differentiable, such that
$|\psi'(t^-)|$ is Riemann integrable with $\psi'(t^-)$ the
left-derivative of $\psi$ in $t$, with support included in $[0,1]$
and $ \int_\Real t^p\psi(t)\,dt=\int_0^1
t^p\psi(t)\,dt=0~~\mbox{for}~~p=0,1.$. Let $A:=\big
|\sup_jH_j^*-\inf_jH_j^*\big |$. If $A<1/2$, under Assumption C, for
all $0<\kappa<\frac 1 {1+4A}-\frac 1 3$, if $a_N=N^{1/3+\kappa}$ and
$v_N=N^{2/3(1-2A)-\kappa(2+4A)}$ then (\ref{conv_tau}),
(\ref{conv_theta}), (\ref{estim2}) and (\ref{test}) hold.
\end{cor}
\begin{remark}
The dependence of this result on $A$ can be explained by the fact that
$2(\sup_j\alpha_j^*-\inf_j\alpha_j^*)+1$, with
$\alpha_j^*=2H_j^*+1$, has to be smaller than $3$ since $a_N \cdot
N^{-1/3}\limiteNN \infty$. However, Corollary \ref{cor_FBM} is quite surprising:
the smaller $A$, {\it i.e.} the smaller the differences between the
parameters $H_j$, the faster the convergence rates of estimators
$\widehat \tau_j$ to $\tau_j^*$. If the difference between two
successive parameters $H_j$ is too large, the estimators
${\widehat\tau_j}$ do not seem to converge.
This is attributable to the influence of the other segments that is
even deeper than the involved exponents are different (simulations
exhibit this paroxysm in \cite{bardet3}).\\ \vspace{-2mm}
\end{remark}
Thus, the rate of convergence of $\underline{\widehat \tau}$ to
$\underline{\tau}^*$ (in probability) can be $N^{2/3(1-2A)-\kappa'}$
for $0<\kappa'$ as small as one wants when
$a_N=N^{1/3+\kappa'/(2+4A)}$. Results of simulations can be observed
in Table \ref{Table} and Figure \ref{Figure} in a case where
$A=0.3<1/2$.





\end{document}